\documentclass[10pt]{amsart}
\usepackage{amssymb,amscd}
\usepackage{ fullpage}

\numberwithin{equation}{section}

\newtheorem{theorem}{Theorem}[section]

\newtheorem{lemma}[theorem]{Lemma}
\theoremstyle{definition}

\theoremstyle{remark}

\newcommand\cA{\mathcal{A}}

\newcommand\cG{\mathcal{G}}

\newcommand{\R}{\mathbb{R}}
\newcommand{\C}{\mathbb{C}}

\newcommand{\bP} {\mathbb{P}}
\newcommand\lie[1]{\mathfrak{#1}}

\newcommand{\fg}{\lie{g}}

\newcommand{\fk}{\lie{k}}

\def    \inv    {^{-1}}

\newcommand{\Hom} {\operatorname{Hom}}

\begin{document}

\title{Gradient flow of the norm squared of a moment 
map}

\author{ E. Lerman }
\thanks{Supported in
part by DMS-0204448}

\address{University of Illinois, Urbana, IL 61801, USA}
\address{Australian National University, Canberra, ACT 0200, Australia 
(till August 2005)}
\email{lerman@math.uiuc.edu}


\begin{abstract}
We present a proof due to Duistermaat that the gradient flow of the
norm squared of the moment map defines a deformation retract of the
appropriate piece of the manifold onto the zero level set of the
moment map.  Duistermaat's proof is an adaptation of Lojasiewicz's
argument for analytic functions to functions which are locally
analytic.
\end{abstract}

\maketitle
\section{Introduction}

Recall that a 2-form $\sigma$ on a manifold $M$ is \emph{symplectic}
if it is closed and nondegenerate. Recall also that if  a Lie group
$G$ acts on a manifold $M$ then for any element $X$ in the Lie algebra
$\fg$ of $G$ we get a vector field $X_M$ on $M$ defined by:
\[
X _M(m) = \left.\frac{d}{dt}\right|_{t=0} \exp tX \cdot m
\]
for all $m\in M$.  Here $\exp : \fg \to G$ denotes the exponential map
and ``$\cdot$'' the action.  Given an action of a Lie group $G$ on a
manifold $M$ preserving a symplectic form $\sigma$ there may
exist a \emph{moment map} $\mu: M \to \fg^*$ associated with this
action.  It takes its values in $\fg^*:=\Hom (\fg, \R)$, 
in the dual of the Lie algebra of $G$.  The moment map is defined as
follows: for any $X\in \fg$, $\langle \mu, X\rangle$ is a real-valued
function on $M$ ($\langle \cdot, \cdot \rangle$ denotes the canonical
pairing $\fg^* \times \fg \to \R$).  We require that
\[
d \langle \mu , X \rangle = \sigma (X_M, \cdot ) 
\]
for any $X\in \fg$ (the opposite sign convention is also used by many
authors: $ d \langle \mu , X \rangle = -\sigma (X_M, \cdot ) $).
Additionally one often requires that the moment map $\mu $ is
\emph{equivariant}, that is, that it intertwines
the given action of $G$ on $M$ with the coadjoint action of $G$ on $\fg^*$:
\[
\mu (g\cdot m) = g\cdot \mu (m)
\]
for all $g\in G$ and $m\in M$.  An equivariant moment map is unique up
to an addition of a vector in $(\fg^*)^G$, the $G$-fixed subspace of
$\fg^*$.  For example if $G$ is a subgroup of the unitary group
$U(n)$, then its natural action on the projective space $\C P^{n-1}$
not only preserves the imaginary part $\sigma$ of the Fubini-Study
metric but 
is, in fact, Hamiltonian. The associated moment map $\mu$ is given by
\[
\langle \mu ([v]) , X\rangle = \frac{1}{||v||^2} (Xv, v),
\]
where $[v]\in \C P^{n-1}$ denotes the line through a vector $v\in \C^n
\smallsetminus 0$, $(\cdot, \cdot)$ denotes the standard Hermitian inner
product on $\C^n$ and $Xv$ denotes the image of $v$ under $X\in
\fg\subset \Hom (\C^n, \C^n)$ (see for
instance \cite[p.~24]{Kbook}).  If $V\subset \C P^{n-1}$ is a
$G$-invariant smooth subvariety, then the action of $G$ on $V$ also
has a moment map and it is simply the restriction $\mu|_V$ .

Equivariant moment maps can be used to construct new symplectic
manifolds: Suppose that the Lie group $G$ acts freely and properly on
the zero level set $\mu\inv (0)$ of an associated moment map $\mu:M\to
\fg^*$. Then by a theorem of Meyer \cite{Meyer} and, independently,
of Marsden and Weinstein \cite{MW}, 
\[
M/\!/_0 G:= \mu\inv (0)/G
\]
is a symplectic manifold, called the \emph{symplectic quotient} at 0.  If
the action of $G$ on $\mu \inv (0)$ is not free then $M/\!/_0 G$ is a
stratified space with symplectic strata \cite{SL}. 

In many cases one knows the topology of $M$ and one is interested in
understanding the topology of the symplectic quotient $M/\!/_0 G$,
which may be much more complicated.  For instance a projective toric
manifold is a symplectic quotient of a symplectic vector space by a
theorem of Delzant.  Such a manifold has an interesting
cohomology ring, while the vector space one constructs it from is
contractible.

A more interesting example is due to Atiyah and Bott \cite{AB}.
Consider the space $\cA$ of connections on a vector bundle $E$ over a
Riemann surface $\Sigma$.  Then $\cA$ is an affine infinite
dimensional space with a constant coefficients symplectic
form. Moreover the action of the gauge group $\cG$ on $\cA$ is
Hamiltonian and the associated moment map $\mu$ assigns to a
connection $A$ its curvature $F_A$.  Therefore the symplectic quotient
$\cA/\!/_0 \cG$ is (formally) the moduli space of flat connections.
Atiyah and Bott also pointed out that the norm squared $||\mu||^2$ of
the moment map behaves somewhat like a very nice Morse-Bott function.
While the connected components of the critical set of $||\mu||^2$ are
not manifolds, they do have indices.  Moreover, for each critical
component $C$ the set $S_C$ given by
\begin{equation}\label{sc}
S_C := \{ x\in M \mid \omega\text{-limit of the  trajectory
of }-\nabla ||\mu||^2  \text{ is in } C\}
\end{equation}
is a manifold --- the ``stable manifold'' of $C$.  The manifolds
$S_C$ give rise to a decomposition of $\cA$. This was later made
rigorous by Donaldson \cite{Donald} and Daskalopoulos \cite{Das}.  In
the finite dimensional setting these ideas were developed
 by Kirwan \cite{Kbook} and Ness \cite{Ness} (independently of each other).  

The work of Kirwan and Ness had an additional motivation that comes
from the Geometric Invariant Theory (GIT) of Mumford.  Roughly
speaking given a complex projective variety $M\subset \bP (V)$ and a
complex reductive group $G^\C \subset \text{PGL} (V)$ one can form a
new projective variety $M/\!/ G^C$, the GIT quotient of $M$.  It is
obtained by taking a Zariski dense subset $M_{ss}\subset M$ of
``semi-stable points'' and dividing out by $G^\C$:
\[
 M/\!/ G^C := M_{ss}/G^C .
\]
It turns out that the action of the maximal compact subgroup $G$ of
$G^C$ on $M$ is Hamiltonian 
and that the two quotients are equal:
\begin{equation}\label{sta}
 M/\!/ G^C = M/\!/_0 G .
\end{equation}
As far as credit for this observation goes, let me
quote from a wonderful paper by R.\ Bott \cite[p.~112]{Bott}
\begin{quote}
In fact, it is quite distressing to see how long it is taking us
collectively to truly sort out symplectic geometry.  I became aware
of this especially when one fine afternoon in 1980, Michael Atiyah and
I were trying to work in my office at Harvard.  I say trying, because
the noise in the neighboring office made by Sternberg and Guillemin
made it difficult.  So we went next door to arrange a truce and in the
process discovered that we were \emph{grosso modo} doing the same
thing.  Later Mumford joined us, and before the afternoon was over we
saw how Mumford's ``stability theory'' fitted with Morse theory.
\end{quote}
Both Guillemin and Sternberg in \cite{GS-quant} and Ness in
\cite{Ness} credit Mumford for (\ref{sta}).

As I mentioned earlier, it is of some interest in symplectic and
algebraic geometry to understand the topology of symplectic and GIT
quotients.  In particular it is natural to ask whether the stable
manifolds $S_C$ retract onto the corresponding critical sets $C$ under
the flow of $-\nabla ||\mu||^2$.  This is not entirely obvious in
light of the fact that there are functions $f$ on $\R^2$ whose
gradient flows have nontrivial $\omega$-limit sets.  As the referee
pointed out, an example can by found on p.~14 of \cite{PdM}: The
function $f$ in question is given in polar coordinates by
\[
f(r,\theta) =\left\{ 
\begin{array}{lr}
e^{\frac{1}{r^2-1}},  & \text{if }
 r <1;\\
0, & \text{if } r = 1;\\
e^{-\frac{1}{r^2-1}} \sin (1/(r-1) -\theta),  & \text{if } r>0 .
\end{array}
\right.
\]
The whole unit circle $\{ r=1\}$ is the $\omega$-limit set of a
gradient trajectory of $f$.  For this function the gradient flow does
not give rise to a map from $S_C$ to $C$, let alone a retraction.

If the function $f$ is analytic then the $\omega$-limit sets of the
flow $\nabla f$ are single points, as was proved by Lojasiewicz
\cite{Lo2}.  Pushing this idea a bit further and using the results of
Kempf and Ness \cite{KN} Neeman proved that in algebraic setting the
flow of $-\nabla ||\mu||^2$ defines a retraction of the set of
semi-stable points $M_{ss}$ onto the zero set of the moment map
$\mu\inv (0)$ \cite{Neeman}.  See also Schwarz
\cite{Schwarz} for a nice survey.  Note that the moment map itself
does not appear explicitly in these two papers.  The connection is
explained in a paper of Linda Ness \cite{Ness} quoted earlier.  In
the setting of connections of vector bundles over Riemann surfaces
Daskolopoulus
\cite{Das} showed that the gradient flow of the $||\mu||^2$ defines a
continuous deformation retract of the Atiyah-Bott strata $S_C$ onto
the components of the critical set of $||\mu||^2$.  This leaves us
with the question: is it true that for an arbitrary moment map $\mu$
the stable manifolds $S_C$ defined by (\ref{sc}) retract onto the
critical sets $C$ under the flow of $-\nabla ||\mu||^2$?.  The answer,
under reasonable assumptions on $\mu$, is yes.  It has been known to
experts for some time.  It was proved by Duistermaat in the 1980's.
The existence of Duistermaat's proof is even footnoted in \cite{MFK}
on p.~166 (this was kindly pointed out by the referee).  The result
has been used by a number of authors, but until Woodward  wrote
it up in
\cite[Appendix B]{Wood}, there was no widely available proof.  In this
paper I do my best to write out Duistermaat's proof in detail.  I do
not claim any originality.  The proof closely follows ideas of
Lojasiewicz
\cite{Lo2} on the properties of gradient flows of analytic functions
(see \cite{KMP}, p.\ 763 and p.\ 765 for a nice summary).  The norm
squared of a moment map is not necessarily an analytic function, but
it is one locally.  The latter is enough to prove that the Lojasiewicz
gradient inequality (see Lemma~\ref{lem2} below) holds for the norm
squared of the moment map and make the rest Lojasiewicz's argument
work.  Woodward proves the Lojasiewicz inequality directly using the
local normal form theorem for moment maps.   He gets more precise
constants for the rate of convergence of the gradient flow but his
exposition is a bit  terse.\\

We now recapitulate our notation and make things a bit more precise.
Let $(M, \sigma)$ be a connected symplectic manifold with a
Hamiltonian action of a compact Lie group $K$ and an associated
equivariant moment map $\mu: M\to \fk^*$.  We assume throughout the
rest of paper that the moment map is {\em proper}.  We fix an
invariant inner product on $\fk^*$ and a Riemannian metric on $M$
compatible with $\sigma$.
Let
\[
 f (x) = ||\mu (x)||^2 ,
\]
denote the norm squared of the moment map. It is a proper map. We
denote the flow of $-\nabla f$  by
$\phi_t$.  Thus, for any function $h$ and any $x\in M$ we have
\[
\frac{d}{dt} h (\phi_t (x)) = -\nabla f (h) \, (\phi _t (x))  
= -(\nabla f \cdot \nabla h )\, (\phi _t (x)).
\]
Since $f$ is proper, the flow exists of all $t\geq 0$.
Kirwan proved \cite[Theorem~4.16, p. 56]{Kbook} that the set of the
critical points of the function $f$ is a  disjoint
union of path connected closed subsets 
on each of which $f$ takes constant value.  Moreover for each
such component $C$ the corresponding stable set $S_C$ defined by 
\[
S_C := \{ x\in M \mid \omega\text{-limit } \phi_t (x)  \subset C\}
\]
is a smooth manifold.   The main result of the paper is

\begin{theorem}[Duistermaat]\label{thm-main}
Let $(M, \sigma)$, $f$, $C$, $\phi_t$ and $S_C$ be as above.
Then
\begin{enumerate}
\item for each $x\in M$ the $\omega$-limit set of the trajectory 
$ \phi_t (x)$ is a
single point, which we denote by $\phi_\infty (x)$;

\item for each connected component $C$ of critical points of $f$
 the map
\[
\phi: [0, \infty]\times S_C \to C, \quad (t, x) \mapsto \phi_t (x)
\] is a deformation retraction.
\end{enumerate}
\end{theorem}

\subsection*{Acknowledgments} I am grateful to the referee and to Anton
Alekseev for helpful comments.  I thank Amnon Neeman for reading an
early version of the manuscript.  I am also very grateful to Hans
Duistermaat for making his 1980's manuscript available to me back in 1988.

\section{Proof of Duistermaat's theorem}

We prove Theorem~\ref{thm-main} in a number of steps.  We first argue
that moment maps for compact Lie groups are locally real analytic.
This local analyticity is a consequence of the local normal form
theorem of Marle
\cite{Marle} and of Guillemin-Sternberg \cite{GSnormal} for moment
maps and of the fact that compact Lie groups are real analytic.
 To state the local normal form theorem we need to set up some notation.
Let $x\in M$ be a point, $G_x$ be its isotropy group with Lie algebra
$\fg_x$, $G_\alpha$ the isotropy group of $\alpha = \mu (x)$ with Lie
algebra $\fg_\alpha$.  Since the action of $G$ is proper, the isotropy
group $G_x$ is compact.  We can then choose a $G_x$-equivariant splitting
\[
\fg^* = \fg_x^* \times (\fg_\alpha/\fg_x)^* \times (\fg/\fg_\alpha)^*
\]
and thereby the embeddings $\fg_x^* \hookrightarrow \fg^*$ and $
(\fg_\alpha/\fg_x)^*\hookrightarrow \fg^*$.

\begin{theorem}[Marle, Guillemin \& Sternberg]
There is a finite dimensional symplectic representation $V$ of $G_x$,
the associated quadratic homogeneous moment map $\mu_V :V\to \fg_x^*$,
a neighborhood $U$ of the orbit $G\cdot x \subset M$, a neighborhood
$U_0$ of the zero section of the vector bundle 
$G\times _{G_x} \left( (\fg_\alpha/\fg_x)^* \times V\right) \to G/G_x$
and a diffeomorphism $\phi: U_0 \to U$ so that
\[
\mu \circ \phi ([g, \eta, v]) = Ad^\dagger (g) (\alpha + \eta + \mu_V (v))
\]
for all $[g, \eta, v] \in U_0$.  Here $Ad^\dagger$ denotes the
co-adjoint action and $[g, \eta, v]$ denotes the orbit of $(g, \eta
,v) \in G\times \left( (\fg_\alpha/\fg_x)^* \times V\right) $ in the
associated bundle $G\times _{G_x} \left( (\fg_\alpha/\fg_x)^* \times
V\right) $.  
\end{theorem}
Now, since $G$ is a compact Lie group, it is real analytic. A choice
of a local analytic section of $G\to G/G_x$ gives coordinates on
$G\times _{G_x} \left( (\fg_\alpha/\fg_x)^* \times V\right) $ making
$\mu \circ \phi$ into a real  analytic map.

Next we recall the gradient inequality of Lojasiewicz.
\begin{lemma}[Lojasiewicz gradient inequality ] \label{lem2}
If $f$ be a real analytic function on an open set $W\subset R^n$ then
for every critical point $x$ of $f$ there is a neighborhood $U_x$ of
$x$ and constants $c_x> 0$ and $\alpha _x$, $0 < \alpha_x <1$, such
that
\begin{equation}\label{eq**}
||\nabla f (y)|| \geq c_x | f(y) - f(x)|^{\alpha _x}
\end{equation}
for all $y\in U_x$.  Here $||\cdot ||$ denotes the standard Euclidean norm.
\end{lemma}
\begin{proof} This is  Proposition~1 on  p.~92 of \cite{Lo-Ensemble}. 
 Alternatively 
see Proposition~6.8 in \cite{BM}.
\end{proof}

Since any Riemannian metric on a relatively compact subset of $\R^n$
is equivalent to the Euclidean metric, the inequality (\ref{eq**})
holds for an arbitrary Riemannian metric on $\R^n$ with the {\em same}
exponent $\alpha_x$ and possibly different constant $c_x$.

Since the connected component $C$ of the set of critical points of $f$
is compact, it can be covered by finitely many open sets $U_i$ on
which the inequality (\ref{eq**}) holds (with the constants $c_i$ and
$\alpha_i$ depending on the $U_i$). Let $\alpha = \max \alpha_i$ and
$b = f(C)$.  Then for any $y\in U_i$ with $|f(y) - b| < 1$ we have
\[
|f(y) - b|^{\alpha_i} \geq |f(y) - b|^{\alpha}.
\] 
Let $c = \min c_i$ and let $U = \bigcup U_i \cap \{ z\in M \mid |f(z)
- b| < 1\}$.  Then for any $y\in U$ we have $ ||\nabla f (y)||  \geq c|
f(y) - b)|^{\alpha} $ or, equivalently,
\begin{equation}\label{eq****}
||\nabla f (y) ||\cdot | f(y) - b)|^{- \alpha} \geq c.
\end{equation}
By definition of $S_C$,
\[
\lim _{t\to +\infty} f (\phi_t (y)) = b.
\]
for any $y\in S_C$.  
Since $f$ is proper, there is $D>0$ so that 
\[
|f(y) - b| < D \Rightarrow y\in U.
\]
Hence for any $y\in S_C$ there is $\tau =  \tau (y)$ so that 
\[
t\geq \tau (y)  \Rightarrow \phi_t (y) \in U .
\]
Fix $y\in S_C$.  For $t> \tau(y)$
\[
\begin{split}
- \frac{d}{dt} ( f(\phi _t (y)) -b)^{1 -\alpha}  &=
 - (1-\alpha)\, ( f(\phi _t (y)) -b)^{ -\alpha} \,
 (-\nabla f ( f)) (\phi_t (y))\\
&= (1-\alpha)\, ( f(\phi _t (y)) -b)^{ -\alpha} \,||\nabla f(\phi_t
(y))||^2 \\
&\geq (1-\alpha)\, c \,||\nabla f (\phi_t (y))||,
\end{split}
\]
where the last inequality follows by (\ref{eq****}).
Hence for any $t_1 > t_0 > \tau (y)$
\[
\begin{split}
( f(\phi _{t_0} (y)) -b)^{1 -\alpha} -( f(\phi _{t_1} (y)) -b)^{1 -\alpha} &=
 - \int_{t_0}^{t_1}  \frac{d}{dt} ( f(\phi _t (y)) -b)^{1 -\alpha} \, dt \\
 &\geq (1-\alpha ) c \int_{t_0}^{t_1} ||\nabla f(\phi_t (y))|| \, dt 
\end{split}
\]
by the previous inequality.
Setting $c' = \frac{1}{(1-\alpha)c}$ we get:
\begin{lemma}\label{prop}
There are constants $c'>0$ and $0< \alpha <1$ such that for any $t_0 <
t_1$ sufficiently large and any $x\in S_C$
\begin{equation} \label{loj-ineq}
c'\left((f(\phi_{t_0}(x) ) - b)^{1-\alpha} - (f(\phi_{t_1}(x)) -
b)^{1-\alpha} \right)
\geq  \int_{t_0}^{t_1} ||\nabla f (\phi _t (x))|| \, dt
\end{equation}
\end{lemma}

It is now easy to show that the $\omega$-limit of $\phi_t (y)$ as
$t\to +\infty$ is a single point.  Denote by $d$ the distance on $M$
defined by the Riemannian metric.  Then
\begin{equation}\label{eq***}
\begin{split}
d (\phi_{t_0} (y), \phi_{t_1} (y))  
&\leq \int_{t_0}^{t_1} ||\frac{d}{dt} (\phi _t (y))||\, dt\\
&= \int_{t_0}^{t_1} ||\nabla f(\phi _t (y))||\, dt \\
&  \leq  c'\,\, \left(
( f(\phi _{t_0} (y)) -b)^{1 -\alpha} - ( f(\phi _{t_1} (y)) -b)^{1 -\alpha}
\right)
\end{split}
\end{equation}
As $t_0, t_1\to + \infty$ the last expression converges to 0.
Therefore, by the Cauchy criterion, $\lim_{t\to + \infty} \phi_t (y)$
does exist, i.e., the map
\[
\phi_\infty: S_C \to C \quad \phi_\infty (y) : = \lim_{t\to + \infty} \phi_t (y)
\]
is a well-defined map.\\

We argue that $\phi_\infty :S_C \to C$ is continuous.  Given
$x\in S_C$ and $\varepsilon >0$ we want to find $\delta$ such that
\[
d (x, y) < \delta \Rightarrow d( \phi_\infty (x), \phi_\infty (y))
 < \varepsilon 
\]
for any $y\in S_C$.  If we take the limit of both sides of
(\ref{eq***}) as $t_1 \to +\infty$ we get 
\begin{equation} \label{two-star}
d (\phi _t (y), \phi _\infty (y)) \leq c' (f(\phi_t (y)) - b) ^{1 -\alpha}
\end{equation}
for all $y\in S_C$ and all $t$ sufficiently large.
Choose $t> 0$ so that 
\begin{equation}\label{eq2.45}
c' (f(\phi _t (x)) -b )^{1-\alpha} < \varepsilon/4 .
\end{equation}
Then
\begin{equation}\label{a}
d (\phi_t (x), \phi_\infty (x)) <\varepsilon/4 . 
\end{equation}
With $t$ fixed as above, choose $\delta >0$ so that ( $y\in S_C$ and
$d(x, y) <\delta$ ) imply two inequalities:
\begin{equation}\label{b1}
d (\phi_t (x), \phi_t (y)) < \varepsilon/4.
\end{equation}
and 
\begin{equation}\label{b2}
{c'} \left|(f(\phi_t (x)) - b) ^{1 -\alpha} -  
(f(\phi_t (y)) - b) ^{1 -\alpha} \right| < \varepsilon/4 .
\end{equation}
This can be done because both maps $z\mapsto \phi_t (z)$ and $z\mapsto
(f(\phi_t (z)) - b) ^{1 -\alpha}$ are continuous.
Equations (\ref{eq2.45}) and (\ref{b2}) imply that
\begin{equation}\label{b3}
{c'} (f(\phi_t (y)) - b) ^{1 -\alpha} < \varepsilon/2,
\end{equation}
and so, by (\ref{two-star}),
\begin{equation}\label{b4}
d (\phi_\infty (y),\phi_t (y)) < \varepsilon/2.
\end{equation}
Putting (\ref{a}), (\ref{b1}) and (\ref{b3}) together, we get that for
$y\in S_C$, $d(x,y)
< \delta $ $\Rightarrow$
\[
d (\phi_\infty (x), \phi_\infty (y)) \leq d (\phi_\infty (x), \phi_t (x)) + 
d (\phi_t (x), \phi_t (y)) + d (\phi_t (y), \phi_\infty (y))  <
\varepsilon/4 + \varepsilon/4 + \varepsilon/2
= \varepsilon .
\]
 This proves that $\phi_\infty: S_C \to
C$ is continuous.\\

Finally it follows from the argument above that for any $y_0\in S_C$
and any $\varepsilon >0$ there are $\delta>0$ and $\tau >0$ so that
\[
t>\tau  \textrm{ and } d(y,y_0) < \delta \Rightarrow  
d(\phi_t (y), \phi_\infty (y_0) ) < \varepsilon
\]
for all $y\in S_C$.  Consequently
\[ 
\phi : [0, \infty] \times S_C \to S_C, \quad (t, y) \mapsto \phi_t (y)
\]
is continuous.  That is, $S_C$ deformation retracts onto $C$.
This concludes the proof of Theorem~\ref{thm-main}. \hfill $\Box$

\end{document}